\documentclass[a4paper,11pt]{article}
\usepackage{graphicx}
\usepackage{epstopdf}
\usepackage{amsthm}
\usepackage[utf8]{inputenc}
\usepackage[english]{babel}
\newtheorem{theorem}{Theorem}[section]

\theoremstyle{definition}

\usepackage{amssymb, amsmath}
\title{New subclass of harmonic univalent functions defined by integral  operator}
\date{}
\begin{document}
\numberwithin{equation}{section}
\maketitle
\date{}
\begin{center}
	\author{{\bf{G. M. Birajdar}}\vspace{.31cm}\\
		School of Mathematics \& Statistics,\\
		Dr. Vishwanath Karad MIT World Peace University,\\
		Pune (M.S) India 411038\\
		Email: gmbirajdar28@gmail.com}\vspace{0.31cm}\\
	\author{{\bf{N. D. Sangle}}\vspace{0.31cm}\\
		Department of Mathematics,\\
		D. Y. Patil College of Engineering \& Technology,\\
		 Kasaba Bawda, Kolhapur\\
		(M.S.) India 416006\\
		Email: navneet\_sangle@rediffmail.com}\vspace{.5cm}\\
	
\end{center}
\vspace{1cm}
\abstract{}
In this paper, we introduce the subclass $SHP^{-m}(\alpha,\beta)$ using integral operator and give sufficient coefficient conditions for normalized harmonic univalent function in the subclass $SHP^{-m}(\alpha,\beta)$.These conditions are also shown to be necessary when the coefficients are negative.\\



\section{Introduction} 
The class $S_H$ investigated by Clunie and Sheil-Small \cite{B1}.They studied geometric subclasses and established some coefficient
bounds. Several researchers have worked on the class $S_H$ and its subclasses.By introducing new subclasses, Silverman \cite{B13}, Silverman and Silvia \cite{B14}, Jahangiri \cite{B7}, Dixit and Porwal \cite{B4} etc. presented a systematic study of harmonic univalent functions.\\
Motivating the research work done by Jahangiri \cite{B2,B3},
Purohit et al.\cite{B6}, Darus and Sangle \cite{B17}, Ravindar et.al \cite{B16}, Bhoosnurmath and Swamy \cite{B5}, Yalcin \cite{B10}, Al-Shaqsi et.al \cite{B11}, Sangle and Birajdar \cite{B18},   Murugusundaramoorthy \cite{B15},  we introduce new subclasses of harmonic mappings using integral operator.\\
Also, we determine extreme points, coefficient estimates of $SHP^{-m}(\alpha,\beta)$ and  $THP^{-m}(\alpha,\beta)$.\\
Suppose A be family of analytic functions defined in unit disk $U$ and 
${A^0}$ be the class of all normalized analytic functions.
Let the functions $h \in A$ be of the form 
\begin{align}
h(z) = z + \sum\limits_{v \ge 2}^\infty  {{a_v}\,{z^v}}
\end{align}
The integral operator $I^{n}$ in \cite{B3}  for the above functions $h$ defined as 
\begin{align*}
{I^0}\,h(z) = h(z).
\end{align*}
\begin{align*}
{I^1}\,h(z) = I(z)=\int\limits_0^z {h(t){t^{ - 1}}dt;}
\end{align*}
\begin{align*}
{I^m}h(z) = I\left( {{I^{m - 1}}h(z)} \right),\,m \in N = \left\{ {1,2,3,...} \right\} 
\end{align*}
\begin{align}
{I^m}h(z) = z + \sum\limits_{k = 2}^\infty  {{[v]^{ - m}}} {a_v}{z^v}
\end{align}
The harmonic functions can be expressed as $f = h + \overline g$ where $h \in {A^0}$ is given by (1.1) and $g \in A$ has power series expansion:
\begin{equation*}
g(z) = \sum\limits_{v \ge 1}^\infty  {{b_v}\,{z^v}} ,
\left| {{b_1}} \right| < 1\
\end{equation*}
Clunie and Sheil-Small \cite{B1} defined function of form $f = h + \overline g$ that are locally univalent, sense-preserving and harmonic in $U$. A sufficient condition for the functions $f$ to be univalent in $U$ is
$\left| {h'(z)} \right| \ge \left| {g'(z)} \right|$ in $U$.\\ 
A function  $f(z)=h+\overline g$ is harmonic starlike \cite{B9} for $\left|z\right|=r<1$,if 
\begin{align*}
\frac{\partial }{{\partial \theta }}\left( {\arg \left( {f(r{e^{i\theta }})} \right)} \right) = {\mathop{\rm Re}\nolimits} \left\{ {\frac{{zh'(z) - \overline {zg'(z)} }}{{h(z) + \overline {g(z)} }}} \right\} > 0.
\end{align*}
The integral operator in \cite{B3} defined for the harmonic functions $f$ by 
\begin{align}
I^mf(z) = I^mh(z) + {\left({- 1} \right)^m}\,\overline {I^m g(z)}
\end{align}
where $I^m$ is defined by (1.2).\\
Now for $0\le \alpha<1 ,m \in N_{0} $ and $z \in U$ ,suppose that  $SHP^{-m}(\alpha,\beta)$ denote the family of harmonic univalent function $f$ of the form 
$f(z)=h(z)+\overline g(z)$ \\
where 
\begin{align}
h(z) = z + \sum\limits_{v \ge 2}^\infty  {{a_v}\,{z^v}},\quad
g(z) = \sum\limits_{v \ge 1}^\infty  {{b_v}\,{z^v}} ,
\left| {{b_1}} \right| < 1\
\end{align}
\begin{align}
{\mathop{\rm Re}\nolimits}\left\{{(1 -\alpha)\frac{{I^mh(z)+{{\left({-1} \right)}^m}\overline {I^mg(z)} }}{z} + \alpha {{\left[ {I^mh(z) + {{\left({-1}\right)}^m}\overline {I^mg(z)}}\right]}^{'}}}\right\}\ge\beta
\end{align}

We further denote by $THP^{-m}(\alpha,\beta)$ subclass of $SHP^{-m}(\alpha,\beta)$ consisting harmonic functions  $f=h+\overline{g} $ in $THP^{-m}(\alpha,\beta)$ so that $h$ and $g$ are of the form 
\begin{align}
h(z)=z-\sum_{v\ge2}^{\infty}{\left|a_{v}\right|}{z^{v}} \  and \   g(z)=\sum_{v\ge1}^{\infty}{\left|b_{v}\right|}{z^{v}}.
\end{align}
In this paper, we will give the sufficient condition for functions $f=h+\overline g$ to be in the subclass $SHP^{-m}(\alpha,\beta)$. It is shown that the coefficient condition is also necessary for the functions in the class $THP^{-m}(\alpha,\beta)$. Coefficient and  distortion bounds, extreme points, convolution conditions and convex combination of this class are obtained. 
\section{Main Results} 
We begin by proving some sharp coefficient inequalities contained in the following theorem.
\begin{theorem}
Let the function $f=h+\overline{g}$ be such that $h$ and $g$ are given by (1.6), furthermore, let 
\begin{align}
\sum\limits_{v \ge 1}^\infty  {\left[v \right]}^{-m}\left( {1 - \alpha  + \alpha v} \right)\left( {\left| {{a_v}} \right| + \left| {{b_v}} \right|} \right) \le 1 - \beta.
\end{align}
Then $f(z)$ is harmonic univalent, sense preserving in $U$ and $ f(z) \in SHP_q^m(\alpha,\beta)$.
\end{theorem} 
\textbf{Proof:}
For $\left| {{z_1}} \right| \le \left| {{z_1}} \right| < 1$, we have by equation (2.1)
\begin{align*}
\left| {f({z_1}) - f({z_2})} \right| &\ge \left| {h({z_1}) - h({z_2})} \right| - \left| {g({z_1}) - g({z_2})} \right|\\
&\ge \left| {{z_1} - {z_2}} \right|\left( {1 - \sum\limits_{v \ge 2}^\infty  {v\left| {{a_v}} \right|{{\left| {{z_2}} \right|}^{v - 1}}}  - \sum\limits_{v \ge 1}^\infty  {v\left| {{b_v}} \right|{{\left| {{z_2}} \right|}^{v - 1}}} } \right)\\
&\ge \left| {{z_1} - {z_2}} \right|\left( {1 - \sum\limits_{v \ge 2}^\infty  {v\left( {\left| {{a_v}} \right| + \left| {{b_v}} \right|} \right){{\left| {{z_2}} \right|}^{v - 1}}}  + \left| {{b_1}} \right|} \right)\\
&\ge \left| {{z_1} - {z_2}} \right|\left( {1 - \sum\limits_{v \ge 2}^\infty  {v\left( {\left| {{a_v}} \right| + \left| {{b_v}} \right|} \right)}  + \left| {{b_1}} \right|} \right)\\
&\ge \left| {{z_1} - {z_2}} \right|\left( {1 - \sum\limits_{v \ge 2}^\infty  {\left[v \right]} ^{-m}\left( {1 - \alpha  + \alpha v} \right)\left( {\left| {{a_v}} \right| + \left| {{b_v}} \right|} \right) + \left| {{b_1}} \right|} \right)\\
&\ge \left| {{z_1} - {z_2}} \right|\left( {1 - (1 - \beta ) - \left| {{b_1}} \right| + \left| {{b_1}} \right|} \right)\\
&\ge \left| {{z_1} - {z_2}} \right|\left| \beta  \right|\ge0. \\
\end{align*}
Hence, $f(z)$ is univalent in $U$. $f(z)$ is sense preserving in $U$. This is because
\begin{align*}
\left| {h'(z)} \right| &\ge 1 - \sum\limits_{v \ge 2}^\infty  {u\left| {{a_v}} \right|{{\left| z \right|}^{v - 1}}}\\
&> 1 - \sum\limits_{v \ge 2}^\infty  {v\left| {{a_v}} \right|} \\
&> 1 - \sum\limits_{v \ge 2}^\infty  {\left[v \right]^{-m}} \left( {1 - \alpha  - \alpha  v} \right)\left| {{a_v}} \right|\\
&\ge \sum\limits_{v \ge 1}^\infty  {\left[v \right]^{-m}} \left( {1 - \alpha  - \alpha v} \right)\left| {{b_u}} \right|\\
&\ge \sum\limits_{v \ge 1}^\infty  {\left[ v \right]^{-m}} \left( {1 - \alpha  - \alpha v} \right)\left| {{b_v}} \right|{\left| z \right|^{v - 1}}\\
&>\sum\limits_{v \ge 1}^\infty  {v\left| {{b_v}} \right|{{\left| z \right|}^{v - 1}}}\\
\left| {h'(z)} \right| &\ge \left| {g'(z)} \right|. 
\end{align*}
Now, we show that $ f(z) \in SHP^m(\alpha,\beta)$.Using the fact that 
${\mathop{\rm Re}\nolimits} \left\{ w \right\} \ge \beta$ if and only if 
$\left| {1 - \beta  + w} \right| \ge \left| {1 + \beta  - w} \right|$,
its sufficient to show that 
\begin{align*}
\left|1-\beta+{\mathop{\rm Re}\nolimits}\left\{{(1 -\alpha)\frac{{I^mh(z)+{{\left({-1} \right)}^m}\overline {I^mg(z)} }}{z} + \alpha {{\left[ {I^mh(z) + {{\left({-1}\right)}^m}\overline {I^mg(z)}}\right]}^{'}}}\right\}\right|
\end{align*}
\begin{align}
-\left|1+\beta-{\mathop{\rm Re}\nolimits}\left\{{(1 -\alpha)\frac{{I^mh(z)+{{\left({-1} \right)}^m}\overline {I^mg(z)} }}{z} - \alpha {{\left[ {I^mh(z) + {{\left({-1}\right)}^m}\overline{I^mg(z)}}\right]}^{'}}}\right\}\right|\ge0.
\end{align}
\begin{align*}
= \left| {2 - \beta  + \sum\limits_{v \ge 2}^\infty  {\left[v \right]^{-m}\left( {1 - \alpha  + \alpha v} \right){a_v}{z^{v - 1}} + {{\left( { - 1} \right)}^m}\sum\limits_{v \ge 1}^\infty  {\left[v \right]^{-m}\left( {1 - \alpha  + \alpha v} \right){b_v}{z^{v - 1}}} } } \right|\\
- \left| {\beta  - \sum\limits_{v \ge 2}^\infty  {\left[v \right]^{-m}\left( {1 - \alpha  + \alpha v} \right){a_u}{z^{v - 1}} - {{\left( { - 1} \right)}^m}\sum\limits_{v \ge 1}^\infty  {\left[v \right]^{-m}\left( {1 - \alpha  + \alpha v} \right){b_v}{z^{v - 1}}} } } \right|
\end{align*}
\begin{align*}
\ge 2\left| {(1 - \beta ) - \sum\limits_{v \ge 2}^\infty  {\left[v \right]^{-m}\left( {1 - \alpha  + \alpha v} \right){a_v}{{\left| z \right|}^{v - 1}} + \sum\limits_{v \ge 1}^\infty  {\left[v \right]^{-m}\left( {1 - \alpha  + \alpha v} \right){b_v}{{\left| z \right|}^{v - 1}}} } } \right|
\end{align*}
\begin{align*}
> 2\left| {(1 - \beta ) - \sum\limits_{v \ge 2}^\infty  {\left[v \right]^{-m}\left( {1 - \alpha  + \alpha v} \right){a_v} + \sum\limits_{v \ge 1}^\infty  {\left[v \right]^{-m}\left( {1 - \alpha  + \alpha v} \right){b_v}} } } \right|.
\end{align*}
This last expression is non-negative by (2.1).\\
The harmonic mappings
\begin{align}
f(z) = z + \sum\limits_{v \ge 2}^\infty  {\frac{{1 - \beta }}{{\left[ v \right]^{-m}\left( {1 - \alpha  + \alpha v} \right)}}{x_v}{z^v} + } \overline {\sum\limits_{v \ge 1}^\infty  {\frac{{1 - \beta }}{{\left[ v \right]^{-m}\left( {1 - \alpha  + \alpha v} \right)}}{y_v}{z^v}} }
\end{align}
where $\sum\limits_{v \ge 2}^\infty  {\left| {{x_v}} \right|}  + \sum\limits_{v \ge 1}^\infty  {\left| {{y_v}} \right| = 1}$ shows that the coefficient bound given by (2.1) is sharp.\\
The functions of the form (2.3) are in $SHP^{-m}(\alpha,\beta)$ because 
\begin{align*}
\sum\limits_{v \ge 1}^\infty  {\left[v \right]^{-m}\left( {1 - \alpha  + \alpha v} \right)(\left| {{a_v}} \right|}  + \left| {{b_v}} \right|) = 1 + (1 - \beta )\sum\limits_{v \ge 2}^\infty  {\left| {{x_v}} \right|}  + \sum\limits_{v \ge 1}^\infty  {\left| {{y_v}} \right| = 2 - \beta }.
\end{align*}
\begin{theorem}
Let the function $f=h+\overline{g}$ be such that $h$ and $g$ are given by (1.6) then $f(z) \in THP^{-m}(\alpha,\beta)$ if and only if 
\begin{align}
\sum\limits_{v \ge 1}^\infty  {\left[v \right]_q^m\left( {1 - \alpha  + \alpha v} \right)(\left| {{a_v}} \right|+\left|b_v\right|)}\le 2 - \beta.
\end{align}
\textbf{Proof:}
The ‘if part’ follows from Theorem 2.1 upon noting that the functions $h(z)$ and $g(z)$ in $f(z) \in SHP^{-m}(\alpha,\beta)$ are of the form (1.6), then  $f(z) \in THP^{-m}(\alpha,\beta)$.\\
For the ‘only if’ part, we show that if  $f(z) \in THP^{-m}(\alpha,\beta)$ then the condition (2.4) holds. Note that a necessary and sufficient condition for $f=h+\overline{g}$ given by (1.6) be in  $THP^{-m}(\alpha,\beta)$ is that
\begin{align*}
{\mathop{\rm Re}\nolimits}\left\{{(1 -\alpha)\frac{{I^mh(z)+{{\left({-1} \right)}^m}\overline {I^mg(z)} }}{z} + \alpha {{\left[ {I^mh(z) + {{\left({-1}\right)}^m}\overline {I^mg(z)}}\right]}^{'}}}\right\}>\beta
\end{align*}
or, equivalently
\begin{align*}
{\mathop{\rm Re}\nolimits} \left\{ {1 - \sum\limits_{v \ge 2}^\infty  {\left[v \right]^{-m}\left( {1 - \alpha  + \alpha u} \right)\left| {{a_v}} \right|{z^{v - 1}} - \sum\limits_{v \ge 1}^\infty  {\left[v \right]^{-m}\left( {1 - \alpha  + \alpha v} \right)\left| {{b_v}} \right|} \,{z^{v - 1}}} } \right\} > \beta.
\end{align*}
If we choose z to be real and $z \to {1^ - }$, we get 
\begin{align*}
1 - \sum\limits_{v \ge 2}^\infty  {\left[v \right]^{-m}\left( {1 - \alpha  + \alpha v} \right)\left| {{a_v}} \right| - \sum\limits_{v \ge 1}^\infty  {\left[v \right]^{-m}\left( {1 - \alpha  + \alpha v} \right)\left| {{b_v}} \right|} \,}  \ge \beta
\end{align*}
this is precisely the assertion of (2.4).
\begin{theorem}
If $f(z) \in THP^{-m}(\alpha,\beta)$ , $\left| z \right| = r < 1$ then
\begin{align*}
\left| {f(z)} \right| \le \left( {1 + \left| {{b_1}} \right|} \right)r + \frac{1}{{\left[ 2 \right]^{-m}\,\left( {1 + \alpha } \right)}}\left( {1 - \left| {{b_1}} \right| - \beta } \right){r^2}
\end{align*}
and
\begin{align}
\left| {f(z)} \right| \ge \left( {1 - \left| {{b_1}} \right|} \right)r - \frac{1}{{\left[ 2 \right]^{-m}\,\left( {1 + \alpha } \right)}}\left( {1 - \left| {{b_1}} \right| - \beta } \right){r^2}
\end{align}
\end{theorem}
\textbf{Proof:}
Taking the absolute values of $f(z)$, we obtain
\begin{align*}
\left| {f(z)} \right| &\le \left( {1 + \left| {{b_1}} \right|} \right)r + \sum\limits_{v \ge 2}^\infty  {\left( {\left| {{a_v}} \right| + \left| {{b_v}} \right|} \right)} \,{r^v}\\
&\le \left( {1 + \left| {{b_1}} \right|} \right)r + \sum\limits_{v \ge 2}^\infty  {\left( {\left| {{a_v}} \right| + \left| {{b_v}} \right|} \right)} \,{r^2}\\
&\le \left( {1 + \left| {{b_1}} \right|} \right)r + \frac{1}{{\left[v \right]^{-m}\left( {1 - \alpha  + \alpha v} \right)}}\sum\limits_{v \ge 2}^\infty  {\left[v \right]^{-m}\left( {1 - \alpha  + \alpha v} \right)\left( {\left| {{a_v}} \right| + \left| {{b_v}} \right|} \right)} \,{r^2}\\
&\le \left( {1 + \left| {{b_1}} \right|} \right)r + \frac{1}{{\left[ 2 \right]^{-m}\left( {1 + \alpha } \right)}}\sum\limits_{v \ge 2}^\infty  {\left[v \right]^{-m}\left( {1 - \alpha  + \alpha v} \right)\left( {\left| {{a_v}} \right| + \left| {{b_v}} \right|} \right)} \,{r^2}\\
&\le \left( {1 + \left| {{b_1}} \right|} \right)r + \frac{1}{{\left[ 2 \right]^{-m}\,\left( {1 + \alpha } \right)}}\left( {1 - \left| {{b_1}} \right| - \beta } \right){r^2}
\end{align*}
and
\begin{align*}
\left| {f(z)} \right| &\ge \left( {1 - \left| {{b_1}} \right|} \right)r - \sum\limits_{v \ge 2}^\infty  {\left( {\left| {{a_v}} \right| + \left| {{b_v}} \right|} \right)} \,{r^v}\\
&\ge \left( {1 - \left| {{b_1}} \right|} \right)r - \sum\limits_{v \ge 2}^\infty  {\left( {\left| {{a_v}} \right| + \left| {{b_v}} \right|} \right)} \,{r^2}\\
&\ge \left( {1 - \left| {{b_1}} \right|} \right)r - \frac{1}{{\left[v \right]^{-m}\left( {1 - \alpha  + \alpha v} \right)}}\sum\limits_{v \ge 2}^\infty  {\left[v \right]^{-m}\left( {1 - \alpha  + \alpha v} \right)\left( {\left| {{a_v}} \right| + \left| {{b_v}} \right|} \right)} \,{r^2}\\
&\ge \left( {1 - \left| {{b_1}} \right|} \right)r - \frac{1}{{\left[ 2 \right]^{-m}\left( {1 + \alpha } \right)}}\sum\limits_{v \ge 2}^\infty  {\left[v \right]^{-m}\left( {1 - \alpha  + \alpha v} \right)\left( {\left| {{a_v}} \right| + \left| {{b_v}} \right|} \right)} \,{r^2}\\
&\ge \left( {1 - \left| {{b_1}} \right|} \right)r - \frac{1}{{\left[ 2 \right]^{-m}\,\left( {1 + \alpha } \right)}}\left( {1 - \left| {{b_1}} \right| - \beta } \right){r^2}.
\end{align*}
For the functions 
\begin{align*}
f(z) = z + \left| {{b_1}} \right|\overline z \, - \frac{1}{{\left[ 2 \right]_q^m\left( {1 + \alpha } \right)}}\left( {1 - \left| {{b_1}} \right| - \beta } \right)\,\,{\overline z ^2}
\end{align*}
and
\begin{align*}
f(z) = z - \left| {{b_1}} \right|z - \frac{1}{{\left[ 2 \right]^{-m}\left( {1 + \alpha } \right)}}\left( {1 - \left| {{b_1}} \right| - \beta } \right)\,\,{z^2}.
\end{align*}
For $\left|b_1\right|\le 1-\beta$ shows that the bounds given in Theorem 2.3 are sharp.
\end{theorem}
Next, we determine the extreme points of the closed convex hulls of $THP^{-m}(\alpha,\beta)$ denoted by $clcoTHP^{-m}(\alpha,\beta)$.
\begin{theorem}
A function $f(z) \in clco\ THP^{-m}(\alpha ,\beta)$ if and only if 
\begin{equation}
f(z) = \sum\limits_{v \ge 1}^\infty  {\left( {{\mu _v}{h_v} + {\eta _v}{g_v}} \right)}
\end{equation}
\begin{equation*}
h_1 (z)=z 
\end{equation*}
where 
\begin{equation*}
{h_v}\left( z \right) = z - \frac{{1 - \beta }}{{\left[v \right]^{-m}\left( {1 - \alpha  + \alpha v} \right)}}{z^v},\,\,(v = 2,3,...) 
\end{equation*}
\begin{equation*}
{g_v}\left( z \right) = z - \frac{{1 - \beta }}{{\left[v \right]^{-m}\left( {1 - \alpha  + \alpha v} \right)}}{\overline z ^v},\,\,(v = 1,2,3,...)
\end{equation*}
\begin{equation*}
\sum\limits_{v \ge 1}^\infty  {\left( {{\mu _v} + {\eta _v}} \right) = 1}, \mu_v\ge0.
\end{equation*}
\end{theorem}
In particular, the extreme points of $THP^{-m}(\alpha,\beta)$ are $\left\{ {{h_v}} \right\}$ and  $\left\{ {{g_v}} \right\}$.\\
\textbf{Proof:}
For the functions $f(z)$ of the form (2.6), we have 
\begin{equation*}
f(z) = \sum\limits_{v \ge 1}^\infty  {\left( {{\mu _v}{h_v} + {\eta _v}{g_v}} \right)}
\end{equation*}
\begin{equation*}
= \sum\limits_{v \ge 1}^\infty  {\left( {{\mu _v} + {\eta _v}} \right)z - } \sum\limits_{v \ge 2}^\infty  {\frac{{1 - \beta }}{{\left[v \right]^{-m}\left( {1 - \alpha  + \alpha v} \right)}}} {\mu _v}{z^v}+{\left( { - 1} \right)^m}\sum\limits_{v \ge 1}^\infty  {\frac{{1 - \beta }}{{\left[v \right]^{-m}\left( {1 - \alpha  + \alpha v} \right)}}} {\eta _v}{\overline z ^v}
\end{equation*}
then
\begin{align*}
\begin{array}{l}
\\
\sum\limits_{v \ge 2}^\infty  {\frac{{\left[ v \right]^{-m}\left( {1 - \alpha  + \alpha v} \right)}}{{1 - \beta }}} \left( {\frac{{1 - \beta }}{{\left[v \right]^{-m}\left( {1 - \alpha  + \alpha v} \right)}}{\mu _v}} \right) + \sum\limits_{v \ge 1}^\infty  {\frac{{\left[v \right]^{-m}\left( {1 - \alpha  + \alpha v} \right)}}{{1 - \beta }}} \left( {\frac{{1 - \beta }}{{\left[v \right]^{-m}\left( {1 - \alpha  + \alpha v} \right)}}{\eta _v}} \right)
\end{array}
\end{align*}
\quad  \quad \quad $= \sum\limits_{v \ge 2}^\infty  {{\mu _v} + } \sum\limits_{v \ge 1}^\infty  {{\eta _v}}  = 1 - {\mu _1} \le 1$\\
and so $f(z) \in clco \  THP^{-m}(\alpha,\beta)$.\\
Conversely ,suppose that $f(z) \in clco\ THP^{-m}(\alpha ,\beta)$.consider \\
\begin{equation*}
{\mu _v} = \frac{{\left[v \right]^{-m}\left( {1 - \alpha  + \alpha v} \right)}}{{1 - \beta }}\left| {{a_v}} \right|\,\,\,,(v = 2,3,...)
\end{equation*}
and
\begin{equation*}
{\eta _v} = \frac{{\left[v \right]^{-m}\left( {1 - \alpha  + \alpha v} \right)}}{{1 - \beta }}\left| {{b_v}} \right|\,\,\,,(v = 1,2,3,...)
\end{equation*}
Then note that by Theorem 2.2 $0 \le {\mu _v} \le 1\,\,(v = 2,3,...)$,
and $0 \le {\eta _v} \le 1\,\,(v =1,2,3,...)$\\
We define ${\mu _1} = 1 - \sum\limits_{v \ge 2}^\infty  {{\mu _v} + } \sum\limits_{v \ge 1}^\infty  {{\eta _v}}$ and note that, by Theorem 2.2 ${\mu _1}\ge0$.\\
Consequently, we obtain \\
$f(z)=\sum\limits_{v \ge 1}^\infty{\left( {{\mu _v}{h_v} + {\eta _v}{g_v}} \right)}$  as required.\\
Using Theorem 2.2, it is easily seen that  $THP^{-m}(\alpha,\beta)$
is convex and closed, so $clcoTHP^{-m}(\alpha,\beta)=THP^{-m}(\alpha,\beta)$.Then the statement of Theorem 2.4 is really for $ f(z) \in THP^{-m}(\alpha,\beta)$.
\begin{theorem}
Each member of  $THP^{-m}(\alpha,\beta)$ maps $U$ on to a starlike domain.
\end{theorem}
\textbf{Proof:}
We only need to show that if $f(z) \in THP^{-m}(\alpha,\beta)$ then 
\begin{align*}
{\mathop{\rm Re}\nolimits} \left\{ {\frac{{zh'(z) - \overline {zg'(z)} }}{{h(z) + \overline {g(z)} }}} \right\} > 0.
\end{align*}
Using the fact that ${\mathop{\rm Re}\nolimits} \left\{ w \right\} > 0$ if and only if $\left| {1 + w} \right| > \left| {1 - w} \right|$,it suffices to show that
\begin{equation*}
\left| {h(z) + \overline {g(z)}  + zh'(z) - \overline {zg'(z)} } \right| - \left| {h(z) + \overline {g(z)}  + zh'(z) + \overline {zg'(z)} } \right|
\end{equation*}
\begin{align*}
&=\left| {2z - \sum\limits_{v \ge 2}^\infty  {(v + 1)\left| {{a_v}} \right|{z^v} + \sum\limits_{v \ge 1}^\infty  {(v - 1)\left| {{b_v}} \right|{z^{ - v}}} } } \right| - \left| {\sum\limits_{v \ge 2}^\infty  {(v - 1)\left| {{a_v}} \right|{z^v} - \sum\limits_{v \ge 1}^\infty  {(v + 1)\left| {{b_v}} \right|{z^{ - v}}} } } \right|\\
&\ge 2\left| z \right| - \left| {\sum\limits_{v \ge 2}^\infty  {(v + 1)\left| {{a_v}} \right|{z^v} + \sum\limits_{v \ge 1}^\infty  {(v - 1)\left| {{b_v}} \right|{z^{ - v}}} } } \right| - \left| {\sum\limits_{v \ge 2}^\infty  {(v - 1)\left| {{a_v}} \right|{z^v} - \sum\limits_{v \ge 1}^\infty  {(v + 1)\left| {{b_v}} \right|{z^{ - v}}} } } \right|\\
&\ge 2\left| z \right|\left\{ {1 - \left( {\sum\limits_{v \ge 2}^\infty  {v\left| {{a_v}} \right|{{\left| z \right|}^{v - 1}} + \sum\limits_{v \ge 1}^\infty  {v\left| {{b_v}} \right|{{\left| z \right|}^{v - 1}}} } } \right)} \right\}\\
&\ge 2\left| z \right|\left\{ {1 - \left( {\sum\limits_{v \ge 2}^\infty  {\left[v \right]^{-m}\left( {1 - \alpha  + \alpha v} \right)\left| {{a_v}} \right| + \sum\limits_{v \ge 1}^\infty  {\left[ v \right]^{-m}\left( {1 - \alpha  + \alpha v} \right)\left| {{b_v}} \right|} } } \right)} \right\}\\
&\ge 2\left| z \right|\left[ {1 - \left( {1 - \beta } \right)} \right]=2\left| z \right|\beta\\
&\ge0.
\end{align*}
\begin{theorem}
If $f(z) \in THP^{-m}(\alpha,\beta)$ then $f(z)$ is convex in the disc 
\begin{equation*}
\left| z \right| < \mathop {\min }\limits_v {\left[ {\frac{{1 - \beta  - \left| {{b_1}} \right|}}{v}} \right]^{\frac{1}{{v - 1}}}},\,\,v = 2,3,...,1 - \beta  > \left| {{b_1}} \right|
\end{equation*}
\end{theorem}
\textbf{Proof:}
Let $f(z) \in THP^{-m}(\alpha,\beta)$ and let $r$ be fixed such that $0<r<1$,then if 
 $ r^{-1}f(rz) \in THP^{-m}(\alpha,\beta)$ and we have 
\begin{align*}
\sum\limits_{v \ge 2}^\infty  {{v^2}\left( {\left| {{a_v}} \right| + \left| {{b_v}} \right|} \right)} {r^{v - 1}} &= \sum\limits_{v \ge 2}^\infty  {v\left( {\left| {{a_v}} \right| + \left| {{b_v}} \right|} \right)} \left( {v{r^{v - 1}}} \right)\\
&\le \sum\limits_{v \ge 2}^\infty  {\left[v \right]^{-m}\left( {1 - \alpha  + \alpha v} \right)\left( {\left| {{a_v}} \right| + \left| {{b_v}} \right|} \right)} \left( {v{r^{v - 1}}} \right)\\
&\le 1 - \beta  - \left| {{b_1}} \right|.
\end{align*}
Provided $v{r^{v - 1}} \le 1 - \beta  - \left| {{b_1}} \right|$,which is true 
\begin{align*}
r < \mathop {\min }\limits_v {\left[ {\frac{{1 - \beta  - \left| {{b_1}} \right|}}{v}} \right]^{\frac{1}{{v - 1}}}},\,\,v = 2,3,...,\ 1 - \beta  > \left| {{b_1}} \right|.
\end{align*}
Following Ruscheweyh \cite{B8}, we call the set
\begin{equation}
{N_\delta }f(z) = \left\{ {G:G(z) = z - \sum\limits_{v \ge 2}^\infty  {\left| {{C_v}} \right|{z^v} - \sum\limits_{v \ge 1}^\infty  {\left| {{D_v}} \right|{z^{ - v}}\,\,\,and\,\,\sum\limits_{v \ge 1}^\infty  {u\left( {\left| {{a_v} - {C_v}} \right| + \left| {{b_v} - {D_v}} \right|} \right) \le \delta } } } } \right\}
\end{equation}
as the $\delta$-neighborhood of $f(z)$. From (2.7) we obtain
\begin{equation}
\sum\limits_{v \ge 1}^\infty  {v\left( {\left| {{a_v} - {C_v}} \right| + \left| {{b_v} - {D_v}} \right|} \right)}  = \left| {{b_1} - {D_1}} \right| + \sum\limits_{u \ge 2}^\infty  {v\left( {\left| {{a_v} - {C_v}} \right| + \left| {{b_v} - {D_v}} \right|} \right)}  \le \delta.
\end{equation}
\begin{theorem}
Let $f(z) \in THP^{-m}(\alpha,\beta)$ and $\delta \le \beta$.If 
$G \in {N_\delta }(f)$, then $G$ is a harmonic starlike function.
\end{theorem}
\textbf{Proof:}
Let $G(z) = z - \sum\limits_{v \ge 2}^\infty  {\left| {{C_v}} \right|{z^v} - } \sum\limits_{v \ge 1}^\infty  {\left| {{D_v}} \right|{z^v}}  \in {N_\delta }f(z)$, we have 
\begin{align*}
\sum\limits_{v \ge 2}^\infty  {v\left( {\left| {{C_v}} \right| + \left| {{D_v}} \right|} \right) + \left| {{D_1}} \right| \le } \sum\limits_{v \ge 2}^\infty  {v\left( {\left| {{a_v} - {C_v}} \right| + \left| {{b_v} - {D_v}} \right|} \right)}  + \sum\limits_{v \ge 2}^\infty  {v\left( {\left| {{a_v}} \right| + \left| {{b_v}} \right|} \right) + \left| {{D_1} - {b_1}} \right| + } \left| {{b_1}} \right|
\end{align*}
\begin{align*}
\le \sum\limits_{v \ge 2}^\infty  {\left[v \right]^{-m}\left( {1 - \alpha  + \alpha v} \right)\left( {\left| {{a_v} - {C_v}} \right| + \left| {{b_v} - {D_v}} \right|} \right) + \left| {{D_1} - {b_1}} \right| + \left| {{b_1}} \right| + } \sum\limits_{u \ge 2}^\infty  {\left[ v \right]^{-m}\left( {1 - \alpha  + \alpha v} \right)\left( {\left| {{a_v}} \right| + \left| {{b_v}} \right|} \right)}
\end{align*}
$\le \delta  + \left| {{b_1}} \right| + \left( {1 - \beta  - \left| {{b_1}} \right|} \right)$\\
$\le 1$.\\
Hence, $G(z)$ is a harmonic starlike function.\\
For next theorem, we require to define the convolution of two harmonic
functions. For harmonic functions of the form
\begin{align*}
f(z) = z - \sum\limits_{v \ge 2}^\infty  {\left| {{a_v}} \right|{z^v} - } \sum\limits_{v \ge 1}^\infty  {\left| {{b_v}} \right|{{\overline z }^v}}
\end{align*}
and 
\begin{align*}
G(z) = z - \sum\limits_{v \ge 2}^\infty  {\left| {{C_v}} \right|{z^v} - } \sum\limits_{v \ge 1}^\infty  {\left| {{D_v}} \right|{{\overline z }^v}}
\end{align*}
we define the convolution of two harmonic functions $f (z)$ and $G(z)$ as
\begin{align*}
\left( {f * G} \right)(z) = f(z) * G(z) = z - \sum\limits_{v \ge 2}^\infty  {\left| {{a_v}} \right|\left| {{C_v}} \right|{z^u} - } \sum\limits_{v \ge 1}^\infty  {\left| {{b_v}} \right|\left| {{D_v}} \right|{{\overline z }^v}}
\end{align*}
Using above definition, we show that the class $THP^{-m}(\alpha,\beta)$ is closed under convolution.
\begin{theorem}
For $0 \le {\alpha _1} \le {\alpha _2}$ ,$0 \le {\beta _1} \le {\beta _2}$
let $f(z) \in THP^{-m}(\alpha_{2},\beta_{2})$ and $G(z) \in THP_q^m(\alpha_{1},\beta_{1})$.Then 
\begin{align*}
(f*G)(z) \in THP^{-m}(\alpha_{2},\beta_{2})\subset THP^{-m}(\alpha_{1},\beta_{1})
\end{align*}
\end{theorem}
\textbf{Proof:}
Let 
\begin{align*}
f(z) = z - \sum\limits_{v \ge 2}^\infty  {\left| {{a_v}} \right|{z^v} - } \sum\limits_{v \ge 1}^\infty  {\left| {{b_v}} \right|{{\overline z }^v}} \quad \in THP^{-m}(\alpha_{2},\beta_{2})\\
G(z) = z - \sum\limits_{v \ge 2}^\infty  {\left| {{C_v}} \right|{z^v} - } \sum\limits_{v \ge 1}^\infty  {\left| {{D_v}} \right|{{\overline z }^u}} \quad \in THP^{-m}(\alpha_{1},\beta_{1}).
\end{align*}
Then the convolution $(f*G)$ is given by (2.9).We wish to show that the coefficient of $(f*G)$ satisfies the required condition given in Theorem 2.2.\\
For  $G(z) \in THP^{-m}(\alpha_{1} ,\beta_{1})$,we note that $\left| {{C_v}} \right| < 1$ and $\left| {{D_v}} \right| < 1$.Now, for the convolution function $f*G$, we obtain
\begin{align*}
\sum\limits_{v \ge 2}^\infty  {\frac{{\left[v \right]^{-m}\left( {1 - {\alpha _1} + v{\alpha _1}} \right)}}{{1 - {\beta _1}}}\left| {{a_v}} \right|\left| {{C_v}} \right|}  + \sum\limits_{v \ge 1}^\infty  {\frac{{\left[v \right]^{-m}\left( {1 - {\alpha _1} + v{\alpha _1}} \right)}}{{1 - {\beta _1}}}\left| {{b_v}} \right|\left| {{D_v}} \right|}\\
\le \sum\limits_{v \ge 2}^\infty  {\frac{{\left[v \right]^{-m}\left( {1 - {\alpha _1} + v{\alpha _1}} \right)}}{{1 - {\beta _1}}}\left| {{a_v}} \right|}  + \sum\limits_{v \ge 1}^\infty  {\frac{{\left[v \right]^{-m}\left( {1 - {\alpha _1} + v{\alpha _1}} \right)}}{{1 - {\beta _1}}}\left| {{b_v}} \right|}\\
\le \sum\limits_{v \ge 2}^\infty  {\frac{{\left[v \right]^{-m}\left( {1 - {\alpha _2} + u{\alpha _2}} \right)}}{{1 - {\beta _2}}}\left| {{a_v}} \right|}  + \sum\limits_{v \ge 1}^\infty  {\frac{{\left[v \right]^{-m}\left( {1 - {\alpha _2} + v{\alpha _2}} \right)}}{{1 - {\beta _2}}}\left| {{b_v}} \right|}\le1.
\end{align*}
Since  $0 \le {\alpha _1} \le {\alpha _2}$, $0 \le {\beta _1} \le {\beta _2}$
let $f(z) \in THP^{-m}(\alpha_{2},\beta_{2})$, thus \\
$(f*G)(z) \in THP^{-m}(\alpha_{2},\beta_{2})\subset THP^{-m}(\alpha_{1} ,\beta_{1})$.

\pagebreak

\end{document}